\documentclass[a4paper,12pt]{article}

\usepackage{amsmath, amsfonts}
\usepackage{amssymb}
\usepackage{amsthm} 
\usepackage[T1]{fontenc}
\usepackage[bf,small,tableposition=top]{caption}
\usepackage{subfig}
\usepackage{setspace}
 \usepackage{multirow}
\usepackage{graphicx}
\usepackage{epstopdf}
\usepackage{array}
\usepackage{lscape}
\usepackage{booktabs,caption,fixltx2e}
\usepackage[flushleft]{threeparttable}
\usepackage{rotating}
\usepackage{natbib}
\usepackage{mathtools}
\usepackage{stmaryrd}
\usepackage{textcomp}

\theoremstyle{plain}
   \newtheorem{thm}{Theorem}

   \newtheorem{cor}{Corollary}
   
\begin{document}

\begin{center}
 \Large\bf{Simultaneous best linear invariant prediction of future order statistics for location-scale and scale families and associated optimality properties}
\end{center}
\begin{center}
\bf {\footnotesize  Narayanaswamy Balakrishnan$^1$ and Ritwik Bhattacharya$^2$}
\end{center}
 \begin{center}
\textit{\scriptsize $^1$Department of Mathematics and Statistics, McMaster University, Hamilton, ON L8S 4K1, Canada}\\
 
\textit{\scriptsize $^2$Department of Industrial Engineering, School of Engineering and Sciences, Tecnol\'{o}gico de Monterrey, Quer\'{e}taro 76130,  M\'{e}xico\\}

 \end{center}

\begin{abstract}
In this article, we first derive an explicit expression for the marginal best linear invariant predictor (BLIP) of an unobserved future order statistic based on a set of early observed ordered statistics. We then derive the joint BLIPs of two future order statistics and prove that the joint predictors are trace-efficient as well as determinant-efficient linear invariant predictors. More generally, the BLIPs are shown to possess complete mean squared predictive error matrix dominance property in the class of all linear invariant predictors of two future unobserved order statistics. Finally, these results are extended to the case of simultaneous BLIPs of any $\ell$ future order statistics. Both scale and location-scale families of distributions are considered as the parent distribution for the development of results. 
\end{abstract}


\section{Introduction} \label{sec1}
\paragraph{}
This article deals with the issues of best linear invariant predictor in order statistics. The mathematical formulation of this problem is as follows. Let us consider a continuous distribution with probability density function
\begin{equation}\label{pdf}
	\frac{1}{\sigma} f\left(\frac{x-\mu}{\sigma} \right),
\end{equation}
where $\mu$ and $\sigma$ are the location and scale parameters, respectively. Suppose the first $r$ order statistics (that is, a Type-II right censored sample) $$X_{1:n}<X_{2:n}<\ldots<X_{r:n},$$ out of a sample of size $n$ from (\ref{pdf}), are observed. Then we are interested in predicting the $(n-r)$ unobserved future order statistics $X_{r+1:n}, X_{r+2:n},\ldots, X_{n:n},$ based on the first $r$ observed order statistics. \\

The first prediction problem was discussed by \cite{Goldberger_1962} within the setup of generalized regression model in which the explicit expression of  the marginal BLUP of unobserved quantity was derived. \cite{Kaminsky_1975} extended the idea in the setup of ordered data. Several interested properties of marginal BLUP was discussed by \cite{Rao_1997}, \cite{Doganaksoy_1997} and \cite{Rao_book1, Bala_2003}. The prediction problem arise in finding the system failure time of a $n$-component parallel system where the first few component failure times are observed. Also, the prediction in order data is applicable in outlier detection \citep[see][]{Balasooriya_1989}. An alternative to the point prediction is the interval prediction which is beyond the scope of this article. However, a brief review on interval prediction can be found in  \cite{Patel_1989}. Coming back to the problems on point prediction, for an extensive review on the point prediction in order statistics, one may refer to \cite{Kaminsky_1998}. The best linear invariant predictor (BLIP) is a larger class of predictor than BLUP through which a reduction in mean squared predictive error is possible. The best linear invariant estimators (BLIEs) of unknown scale and location parameters was treated by \cite{Mann_1969}. BLIP of unobserved order statistics was found by \cite{Kaminsky_1975a}. While all the above referred articles deal with the marginal predictor problems, simultaneous prediction problems (both BLUP and BLIP) were ignored in the literature.  This is the prime motivation of this article. \\

Simultaneous prediction problem in order statistics was first attempted by \cite{Bala_2021}. The explicit expressions of joint BLUPs of two future order statistics, obtained by minimizing the determinant of the variance–covariance matrix of the predictors, were presented and the gain in efficiency over marginal BLUPs was established.  Moreover, the non-existence of joint BLUPs of more than two future order statistics was demonstrated. Later on, the simultaneous prediction of any $\ell$ future order statistics, obtained by minimizing the means squared predictive error matrix of the predictors, was discussed by \cite{Bala_2021a}. It has shown that the simultaneous BLUPs are identical with the corresponding marginal BLUPs. In this article, our aim is to determine the simultaneous BLIPs of any $\ell$ future order statistics based on early observed ordered statistics. To begin with, by minimizing the mean squared predictive error of the predictor, we first derived the expression of marginal BLIP of a future order statistics which is independent from the result by \cite{Kaminsky_1975a}. A comparative study between BLUP and BLIP is carried out to demonstrate the performance. It has found that BLIP always yields less mean squared error than that of BLUP. Then, the joint BLIPs of two future order statistics are presented along with some associated properties. In case of simultaneous BLIP, a practical data-driven guideline is provided in order to choose between BLIP and BLUP. Finally, the simultaneous BLIPs of any $\ell$ future order statistics are derived. It is shown that the simultaneous BLIPs are identical with corresponding marginal BLIPs. All these developments are presented under both scale and location-scale family of distributions as the parent distribution for the underlying variables.\\

The rest of this paper is organized as follows. First, we provide a general background on the known results on linear estimation and prediction problems in Section 2. The marginal BLIP case and its associated properties are presented in Section 3.  Simultaneous prediction of two future order statistics is then developed in Section 4.  The complete  mean squared predictive error matrix dominance property of the joint BLIPs are also demonstrated here. The simultaneous prediction of any $\ell$ future order statistics are discussed in Section 5. BLIP in scale family of distributions is presented in Section 6. Finally, some concluding remarks are made in Section 7.

\section{Basic results on linear estimators and predictors}
\paragraph{}
In this section, we present some basic results on linear prediction of  order statistics that are known in the literature. \\
\subsection{BLUEs and BLUPs}
\paragraph{}
Let $\boldsymbol{X} = (X_{1:n},\ldots, X_{r:n})^{\top}_{r\times 1}$ be the available Type-II right censored data from a location-scale family of distributions in (\ref{pdf}). Let us then use the following notation: $\alpha_i$ for the expected value of the standardized order statistic
\begin{equation*}
	Z_{i:n}=  \frac{X_{i:n}-\mu}{\sigma}, ~i\in\{1,\ldots,r\},
\end{equation*}with $\boldsymbol{\alpha}=(\alpha_1, \ldots, \alpha_r)^{\top}_{r\times1}$ and $\sigma^2 \boldsymbol{\Sigma}$ for the variance-covariance matrix of $\boldsymbol{X}$, where $ \boldsymbol{\Sigma}$ is the $r\times r$ covariance matrix of $Z_{i:n}, i\in\{1,\ldots,r\}$, assumed to be positive definite. In addition, let $\boldsymbol{1}$ denote a column vector  $(1,\ldots,1)^{\top}_{r\times 1}$. Using these notations, we then can write
\begin{equation*}
\boldsymbol{X} = \mu\boldsymbol{1} + \sigma\boldsymbol{Z}
\end{equation*}
and
\begin{equation*}
	\mbox{E}[\boldsymbol{X}] = \mu\boldsymbol{1} + \sigma \boldsymbol{\alpha}.
\end{equation*}Thence, upon minimizing the generalized variance
\begin{equation*}
	(\boldsymbol{X}-\mu\boldsymbol{1}-\sigma \boldsymbol{\alpha})^{\top}\boldsymbol{\Sigma}^{-1}(\boldsymbol{X}-\mu\boldsymbol{1}-\sigma \boldsymbol{\alpha})
\end{equation*}with respect to $\mu$ and $\sigma$, the BLUEs  $(\mu^*, \sigma^*)$ of $(\mu, \sigma)$ are obtained as
\begin{eqnarray*}
	\mu^* &=& \frac{1}{\Delta}\{ (\boldsymbol{\alpha}^{\top} \boldsymbol{\Sigma}^{-1} \boldsymbol{\alpha})(\boldsymbol{1}^{\top} \boldsymbol{\Sigma}^{-1}) -  (\boldsymbol{\alpha}^{\top} \boldsymbol{\Sigma}^{-1} \boldsymbol{1})(\boldsymbol{\alpha}^{\top} \boldsymbol{\Sigma}^{-1})\}\boldsymbol{X},\\
	\sigma^* &=& \frac{1}{\Delta}\{ (\boldsymbol{1}^{\top} \boldsymbol{\Sigma}^{-1} \boldsymbol{1})(\boldsymbol{\alpha}^{\top} \boldsymbol{\Sigma}^{-1}) -  (\boldsymbol{1}^{\top} \boldsymbol{\Sigma}^{-1} \boldsymbol{\alpha})(\boldsymbol{1}^{\top} \boldsymbol{\Sigma}^{-1})\}\boldsymbol{X},
\end{eqnarray*}with
\begin{eqnarray}\nonumber
\mbox{Var}(\mu^*) &=& \frac{\boldsymbol{\alpha}^{\top} \boldsymbol{\Sigma}^{-1} \boldsymbol{\alpha}}{\Delta} \sigma^2,~ \mbox{Var}(\sigma^*) ~=~ \frac{\boldsymbol{1}^{\top} \boldsymbol{\Sigma}^{-1} \boldsymbol{1}}{\Delta}\sigma^2,~ \mbox{Cov}(	\mu^*,	\sigma^*) ~=~ -\frac{\boldsymbol{1}^{\top} \boldsymbol{\Sigma}^{-1} \boldsymbol{\alpha}}{\Delta}\sigma^2
\end{eqnarray}	
and
\begin{eqnarray}		\label{bigdelta}
\Delta &=& (\boldsymbol{1}^{\top} \boldsymbol{\Sigma}^{-1} \boldsymbol{1}) (\boldsymbol{\alpha}^{\top} \boldsymbol{\Sigma}^{-1} \boldsymbol{\alpha}) - (\boldsymbol{1}^{\top} \boldsymbol{\Sigma}^{-1} \boldsymbol{\alpha})^2
\end{eqnarray}being the determinant of the matrix
\begin{equation*}
\begin{bmatrix} \boldsymbol{1}^{\top}\boldsymbol{\Sigma}^{-1}\boldsymbol{1} &  \boldsymbol{1}^{\top}\boldsymbol{\Sigma}^{-1}\boldsymbol{\alpha}\\ \boldsymbol{1}^{\top}\boldsymbol{\Sigma}^{-1}\boldsymbol{\alpha} & \boldsymbol{\alpha}^{\top}\boldsymbol{\Sigma}^{-1}\boldsymbol{\alpha}\end{bmatrix}, 
\end{equation*}which is related to the variance-covariance matrix of the BLUEs $(\hat{\mu}, \hat{\sigma})$; see \cite{Cohen_book}, \cite{David_book} and \cite{Nagaraja_book} for pertinent details. \\

Using the general framework of best linear unbiased prediction developed by \cite{Goldberger_1962} for generalized linear regression model, a point BLUP of $X_{s:t},\quad r<s\leq n,$ was derived by \cite{Kaminsky_1975} as
\begin{equation}\label{MBLUP1}
	\tilde{X}_{s:n} = \mu^* + \sigma^*\alpha_s + \boldsymbol{\omega}^{\top}_s\boldsymbol{\Sigma}^{-1}(\boldsymbol{X}-\mu^*\boldsymbol{1} - \sigma^*\boldsymbol{\alpha}),
\end{equation}where $\boldsymbol{\omega}_s = (\omega_1,\ldots, \omega_r)^{\top}_{r\times 1}$, with $\omega_i = \mbox{Cov}(Z_{i:n}, Z_{s:n})$. Consequently, the mean squared predictive error is given by
\begin{equation}\label{MMSPE1}
	\mbox{MSPE}(\tilde{X}_{s:n}) = \sigma^2 \{\mbox{Var}(X_{s:n}) - \boldsymbol{\omega}^{\top}_s\boldsymbol{\Sigma}^{-1} \boldsymbol{\omega}_s + c_{11}\},
\end{equation}where $c_{11} = \mbox{Var}\{(1- \boldsymbol{\omega}^{\top}_s\boldsymbol{\Sigma}^{-1} \boldsymbol{1})\mu^* + (\alpha_s- \boldsymbol{\omega}^{\top}_s\boldsymbol{\Sigma}^{-1} \boldsymbol{\alpha})\sigma^*\}/\sigma^2$. An alternative expression to (\ref{MBLUP1}) of $\tilde{X}_{s:n}$ has been presented recently by \cite{Bala_2021a} as
\begin{equation}\label{MBLUP2}
	\tilde{X}_{s:n} = \boldsymbol{a}^{\top}\boldsymbol{X},
\end{equation}where
\begin{equation*}
	\boldsymbol{a} = \boldsymbol{\Sigma}^{-1}\boldsymbol{\omega}_s + \frac{1}{\Delta} (V_2A_s - V_3B_s)\boldsymbol{\Sigma}^{-1}\boldsymbol{1} +\frac{1}{\Delta} (V_1B_s - V_3A_s) \boldsymbol{\Sigma}^{-1}\boldsymbol{\alpha},
\end{equation*}with $V_1 = \boldsymbol{1}^{\top}\boldsymbol{\Sigma}^{-1}\boldsymbol{1}$, $V_2 = \boldsymbol{\alpha}^{\top}\boldsymbol{\Sigma}^{-1}\boldsymbol{\alpha}$, $V_3=\boldsymbol{1}^{\top}\boldsymbol{\Sigma}^{-1}\boldsymbol{\alpha}$, $A_s = 1-\boldsymbol{1}^{\top}\boldsymbol{\Sigma}^{-1}\boldsymbol{\omega}_s$, $B_s = \alpha_s-\boldsymbol{\alpha}^{\top}\boldsymbol{\Sigma}^{-1}\boldsymbol{\omega}_s$, and $\Delta$ as defined earlier in (\ref{bigdelta}).\\

Proceeding similarly,  \citet[][Theorem 1]{Bala_2021} also derived explicit expressions for the joint best linear unbiased predictors $(\tilde{X}_{s:n}, \tilde{X}_{t:n})$ of $(X_{s:n}, X_{t:n})$, for $r<s<t\leq n$, as 
\begin{equation*}
		\tilde{X}_{s:n} = \boldsymbol{a}^{\top}\boldsymbol{X}\quad\mbox{and}\quad 	\tilde{X}_{t:n} = \boldsymbol{b}^{\top}\boldsymbol{X},
\end{equation*}where
\begin{equation*}
	\boldsymbol{b} = \boldsymbol{\Sigma}^{-1}\boldsymbol{\omega}_t + \frac{1}{\Delta} (V_2A_t - V_3B_t)\boldsymbol{\Sigma}^{-1}\boldsymbol{1} +\frac{1}{\Delta} (V_1B_t - V_3A_t) \boldsymbol{\Sigma}^{-1}\boldsymbol{\alpha},
\end{equation*}with $A_t = 1-\boldsymbol{1}^{\top}\boldsymbol{\Sigma}^{-1}\boldsymbol{\omega}_t$, $B_t = \alpha_t-\boldsymbol{\alpha}^{\top}\boldsymbol{\Sigma}^{-1}\boldsymbol{\omega}_t$ and $\boldsymbol{a}$ defined exactly as in (\ref{MBLUP2}). Moreover, the joint MSPE matrix has been given by these authors to be
	\begin{equation*}
	\sigma^2 \begin{bmatrix} W_{11} &  W_{12}\\ W_{12} & W_{22}\end{bmatrix},
\end{equation*}where 
\begin{eqnarray*}
	W_{11} &=& \omega_{ss} - \boldsymbol{\omega}^{'}_s\boldsymbol{\Sigma}^{-1}\boldsymbol{\omega}_s + \begin{bmatrix} A_s &  B_s\end{bmatrix}\begin{bmatrix} V_{1} &  V_{3}\\ V_{3} & V_{2}\end{bmatrix}^{-1}\begin{bmatrix} A_s\\ B_s \end{bmatrix},\\
	W_{22} &=& \omega_{tt} - \boldsymbol{\omega}^{'}_t\boldsymbol{\Sigma}^{-1}\boldsymbol{\omega}_t + \begin{bmatrix} A_t &  B_t\end{bmatrix}\begin{bmatrix} V_{1} &  V_{3}\\ V_{3} & V_{2}\end{bmatrix}^{-1}\begin{bmatrix} A_t\\ B_t \end{bmatrix},\\
	W_{12} &=& \omega_{st} - \boldsymbol{\omega}^{'}_s\boldsymbol{\Sigma}^{-1}\boldsymbol{\omega}_t + \begin{bmatrix} A_s &  B_s\end{bmatrix}\begin{bmatrix} V_{1} &  V_{3}\\ V_{3} & V_{2}\end{bmatrix}^{-1}\begin{bmatrix} A_t\\ B_t \end{bmatrix}.
\end{eqnarray*}

Further, explicit expressions for the simultaneous BLUPs of any $\ell$ future order statistics have been given in Theorem 3 of \cite{Bala_2021a}. 

\subsection{BLIEs and BLIPs}
\paragraph{}
The best linear invariant estimators of $(\mu, \sigma)$ were first developed by \cite{Mann_1969}, and upon using these results of \cite{Mann_1969}, an expression for the BLIP $\hat{X}_{s:n}$ of $X_{s:n},\quad r<s\leq n$ was presented by \cite{Kaminsky_1975a} as
\begin{equation}\label{BLIP1}
	\hat{X}_{s:n} = \tilde{X}_{s:n} - \frac{c_{12}}{1+c_{22}}\sigma^*,
\end{equation}where $\sigma^2c_{12} = \mbox{Cov}(\sigma^*, \tilde{X}_{s:n} - \boldsymbol{\omega}^{\top}_s\boldsymbol{\Sigma}^{-1}\boldsymbol{X})$ and $\sigma^2c_{22}=\mbox{Var}(\sigma^*)$. The corresponding mean squared predictive error is given by
\begin{equation*}
	\mbox{MSPE}(\hat{X}_{s:n}) = \mbox{MSPE}(\tilde{X}_{s:n}) - \frac{c^2_{12}}{1 + c_{22}}\sigma^2,
\end{equation*}where $\mbox{MSPE}(\tilde{X}_{s:n}) $ is as given in (\ref{MMSPE1}). In the next section, we first present an alternate expression to (\ref{BLIP1}) for the marginal BLIP of $X_{s:n}$ and its MSPE.

\section{Marginal best linear invariant predictor of an order statistic}
\paragraph{}
In this section, we derive an explicit expression for BLIP of $X_{s:n}$ and then discuss some of its properties.
\begin{thm}
	The best linear invariant predictor $\hat{X}_{s:n}$, obtained by minimizing the mean squared predictive error, is of the form $\hat{X}_{s:n} = \boldsymbol{a}^{\top}\boldsymbol{X}$ in which the coefficient $\boldsymbol{a}=(a_1, \ldots, a_r)^{\top}_{r\times 1}$  is given by
	\begin{equation*}
\boldsymbol{a} = \boldsymbol{\Gamma}^{-1} \boldsymbol{\Delta}_{s},
	\end{equation*}where
\begin{equation*}
	\boldsymbol{\Gamma} = \boldsymbol{\Sigma} + (\boldsymbol{\alpha}+\delta\boldsymbol{1})  (\boldsymbol{\alpha}+\delta\boldsymbol{1})^{\top}
\end{equation*}and
\begin{equation*}
	\boldsymbol{\Delta}_{s} = \boldsymbol{\omega}_s + (\alpha_s + \delta)(\boldsymbol{\alpha} + \delta\boldsymbol{1}),
\end{equation*}with $\delta = \frac{\mu}{\sigma}$.
\end{thm}

\noindent\textbf{Proof:} For deriving the BLIP, let us consider the MSPE of $\hat{X}_{s:n}$ given by 
\begin{eqnarray}\nonumber
	\mbox{MSPE}(\hat{X}_{s:n}) &=& E[(\hat{X}_{s:n} - X_{s:n})^2]\\\label{MSPE}
	                                                     &=& E[( \boldsymbol{a}^{\top}\boldsymbol{X}- X_{s:n})^2].
\end{eqnarray}Upon replacing $\boldsymbol{X}$ and $X_{s:n}$ with their corresponding standardized counterparts $\mu\boldsymbol{1} + \sigma\boldsymbol{Z}$ and $\mu + \sigma Z_{s:n}$, respectively, Eq. (\ref{MSPE}) can be simplified as 
\begin{eqnarray}\nonumber
	\mbox{MSPE}(\hat{X}_{s:n})  &=& \mu^2(\boldsymbol{a}^{\top}\boldsymbol{1} - 1)^2 + 2\mu\sigma(\boldsymbol{a}^{\top}\boldsymbol{1} - 1)(\boldsymbol{a}^{\top}\boldsymbol{\alpha} - \alpha_s)\\\nonumber
	&& ~~~~~+\sigma^2\{ \boldsymbol{a}^{\top}\boldsymbol{\Sigma}\boldsymbol{a} -2\boldsymbol{a}^{\top}\boldsymbol{\omega}_s + \omega_{ss}  +  (\boldsymbol{a}^{\top}\boldsymbol{\alpha} - \alpha_s)^2\}\\\nonumber
  &=& \sigma^2\bigg[  \delta^2(\boldsymbol{a}^{\top}\boldsymbol{1} - 1)^2 +2\delta (\boldsymbol{a}^{\top}\boldsymbol{1} - 1)(\boldsymbol{a}^{\top}\boldsymbol{\alpha} - \alpha_s)\\\label{FMSPE}
		&& ~~~~~~~~~~+\boldsymbol{a}^{\top}\boldsymbol{\Sigma}\boldsymbol{a} -2\boldsymbol{a}^{\top}\boldsymbol{\omega}_s + \omega_{ss}  +  (\boldsymbol{a}^{\top}\boldsymbol{\alpha} - \alpha_s)^2\bigg].
\end{eqnarray}Then, the theorem follows readily when we minimize the  $\mbox{MSPE}(\hat{X}_{s:n})$ in (\ref{FMSPE}) with respect to $\boldsymbol{a}$. \\

\begin{cor}
	For $s<t\leq n$, the BLIP $\hat{X}_{t:n}$ of $X_{t:n}$ can be obtained simply by replacing  $\boldsymbol{\Delta}_{s}$ by $\boldsymbol{\Delta}_{t}$ in Theorem 1.  Thus, the predictors $\hat{X}_{s:n}$ and $\hat{X}_{t:n}$ are indeed trace-efficient invariant predictors by their very construction. 
\end{cor}

\subsection{Interpretation of the quantity $\delta$}
\begin{figure}[h]
	\centering
	\includegraphics[width=10cm,height=12cm,angle=-90,scale=1.15]{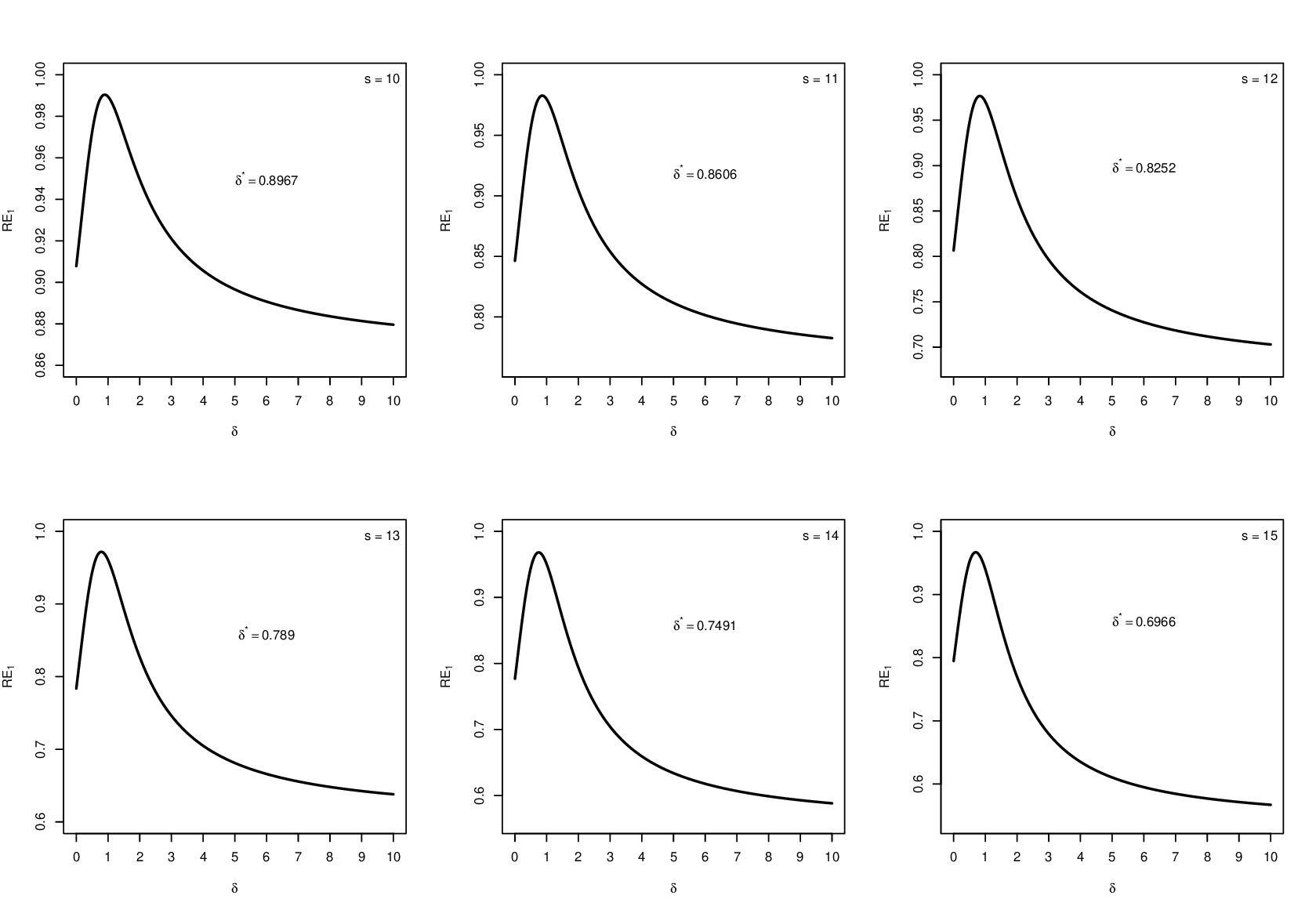}
	\caption{Plot of relative efficiencies of the marginal BLIP against marginal BLUP for the $s$-th order statistic based on $r=9$ and $n=15$. }
	\label{Eff_plot_1}
\end{figure}

\paragraph{}
The quantity $\delta$ plays an important role in the relative performance of the BLIP $\hat{X}_{s:n}$ against the BLUP $\tilde{X}_{s:n}$. Let us define the relative efficiency measure 
\begin{equation*}
	\mbox{RE}_1 = \frac{\mbox{MSPE}(\hat{X}_{s:n}) }{\mbox{MSPE}(\tilde{X}_{s:n})}.
\end{equation*}Note that $\mbox{RE}_1$ is a continuous function of $\delta$. For example, with $n = 15$ and $r = 9$, we have plotted RE$_1$ against $\delta$ in Figure \ref{Eff_plot_1} for $s = 10, 11, 12, 13, 14$ and 15. Note that the range of $\delta$ is $(-\infty, 0)\cup(0, \infty)$. However, due to symmetry, we an focus on the interpretation of the behavior of $\mbox{RE}_1$ in $(0, \infty)$.  Figure \ref{Eff_plot_1} shows the unique maximum of $\mbox{RE}_1$ is attended at some $\delta^{*}$ which can be found numerically. Note that all $\mbox{RE}_1$ are less than 1 indicating that BLIP always possesses less mean squared error than BLUP. When $\delta\rightarrow\delta^*$, $\mbox{RE}_1$ approaches to its maximum value. Therefore, from a practical point of view, one can measure how well BLIP performs better than BLUP by simply estimating $\delta$, say $\hat{\delta}$, using the BLUEs of $\mu$ and $\sigma$ as $\hat{\delta} = \mu^*/\sigma^*$. Thus, for values of $\hat{\delta}$ away from $\delta^{*}$, the BLIP will have better performance meaning less mean squared error.  

\subsection{An illustrative example}
\paragraph{}
Let us consider an environmental lead contamination data 26, 63, 3, 70, 16, 5, 1, 57, 5, 3, 24, 2, 1, 48 and 3 presented by \cite{Bhaumik_2004}. The data were also analyzed by \cite{Krishnamoorthy_2018} and showed that a log-normal distribution fits the data well. In our case, we first take the logarithmic transformation of the data. Then, assuming that the first $r=9$ ordered data are observed, we have $\boldsymbol{X}=(0, 0, 0.693, 1.099, 1.099, 1.099, 1.609,$ $ 1.609, 2.773)^{\top}$. Based on $\boldsymbol{X}$, the BLUEs of $(\mu, \sigma)$ are computed as $(\mu^*, \sigma^*) = (2.253, 1.696)$ which yield $\hat{\delta}=1.328$ for the current data set. We present $\mbox{RE}_1$ for the marginal predictors for $s = 10, 11, 12, 13, 14$ and 15 with $\delta = \hat{\delta}$ in Table 1. The summary of Table 1 indicates that BLIP is better than BLUP when $\hat{\delta}$ is away from $\delta^*$, as expected.
\begin{table}
	\centering
	\caption{Summary of relative efficiencies based on the data presented by \cite{Bhaumik_2004} and \cite{Krishnamoorthy_2018}. Here, $(\mu^*, \sigma^*) = (2.253, 1.696)$ and $\hat{\delta}=1.328$. }
	\begin{tabular}{lllllll}\toprule
	
$s$ & $\delta^*$ & $\hat{X}_{s:n}$ & $\mbox{MSPE}(\hat{X}_{s:n})$ & $\tilde{X}_{s:n}$ & $\mbox{MSPE}(\tilde{X}_{s:n})$ & $\mbox{RE}_1$\\\toprule
10 & 0.8967 & 3.015 & 0.0287 & 3.037 & 0.0293 & 0.9795 \\
11 & 0.8606 & 3.278 & 0.0637 & 3.321 & 0.0664 & 0.9593 \\
12 & 0.8252 & 3.575 & 0.1084 & 3.639 & 0.1157 &  0.9369 \\
13 & 0.7890 & 3.927 & 0.1703 & 4.014 & 0.1855 & 0.9181 \\
14 & 0.7491  & 4.388 & 0.2698& 4.503& 0.3004 & 0.8981 \\
15 & 0.6966 & 5.151  & 0.5037 & 5.305 & 0.5721 & 0.8804\\\bottomrule
	\end{tabular}
\end{table}

\section{Joint best linear invariant predictors of two order statistics}
\paragraph{}
In this section, we derive explicit expressions for the joint best linear invariant predictors of two future order statistics and the corresponding mean squared predictive error matrix.
\begin{thm}
	The joint best linear invariant predictors $\hat{X}_{s:n}$ and $\hat{X}_{t:n}$, for $r<s<t\leq n$, obtained by minimizing the determinant of the mean squared predictive error matrix, are of the form $\hat{X}_{s:n} = \boldsymbol{a}^{\top}\boldsymbol{X}$ and $\hat{X}_{t:n} = \boldsymbol{b}^{\top}\boldsymbol{X}$ in which the coefficients $\boldsymbol{a}=(a_1, \ldots, a_r)^{\top}_{r\times 1}$ and $\boldsymbol{b}=(b_1, \ldots, b_r)^{\top}_{r\times 1}$ are given by	
	\begin{equation*}
		\boldsymbol{a} = \boldsymbol{\Gamma}^{-1} \boldsymbol{\Delta}_{s}\quad\mbox{and}\quad\boldsymbol{b} = \boldsymbol{\Gamma}^{-1} \boldsymbol{\Delta}_{t},
	\end{equation*}where $\boldsymbol{\Gamma} = \boldsymbol{\Sigma} + (\boldsymbol{\alpha}+\delta\boldsymbol{1})  (\boldsymbol{\alpha}+\delta\boldsymbol{1})^{\top}$, $\boldsymbol{\Delta}_{s} = \boldsymbol{\omega}_s + (\alpha_s + \delta)(\boldsymbol{\alpha} + \delta\boldsymbol{1})$ and $\boldsymbol{\Delta}_{t} = \boldsymbol{\omega}_t + (\alpha_t + \delta)(\boldsymbol{\alpha} + \delta\boldsymbol{1})$.
\end{thm}
\noindent\textbf{Proof:} The joint BLIPs will be derived by minimizing the determinant of the mean squared predictive error matrix given by
\begin{equation*}
	 \begin{bmatrix} W_{1} &  W_{3}\\ W_{3} & W_{2}\end{bmatrix},
\end{equation*}where 
\begin{eqnarray*}
W_1 &=& E[(\hat{X}_{s:n} - X_{s:n})^2]\\
         &=& \sigma^2\bigg[  \delta^2(\boldsymbol{a}^{\top}\boldsymbol{1} - 1)^2 +2\delta (\boldsymbol{a}^{\top}\boldsymbol{1} - 1)(\boldsymbol{a}^{\top}\boldsymbol{\alpha} - \alpha_s)\\
         && ~~~~~~~~~~+\boldsymbol{a}^{\top}\boldsymbol{\Sigma}\boldsymbol{a} -2\boldsymbol{a}^{\top}\boldsymbol{\omega}_s + \omega_{ss}  +  (\boldsymbol{a}^{\top}\boldsymbol{\alpha} - \alpha_s)^2\bigg]\\
W_2 &=& E[(\hat{X}_{t:n} - X_{t:n})^2]\\
&=& \sigma^2\bigg[  \delta^2(\boldsymbol{b}^{\top}\boldsymbol{1} - 1)^2 +2\delta (\boldsymbol{b}^{\top}\boldsymbol{1} - 1)(\boldsymbol{b}^{\top}\boldsymbol{\alpha} - \alpha_t)\\
&& ~~~~~~~~~~+\boldsymbol{b}^{\top}\boldsymbol{\Sigma}\boldsymbol{b} -2\boldsymbol{b}^{\top}\boldsymbol{\omega}_t + \omega_{tt}  +  (\boldsymbol{b}^{\top}\boldsymbol{\alpha} - \alpha_t)^2\bigg]\\
W_3 &=& E[(\hat{X}_{s:n} - X_{s:n})(\hat{X}_{t:n} - X_{t:n})]\\   
         &=&   \sigma^2\bigg[   \delta^2(\boldsymbol{a}^{\top}\boldsymbol{1} - 1)(\boldsymbol{b}^{\top}\boldsymbol{1} - 1)+\delta(\boldsymbol{a}^{\top}\boldsymbol{1} - 1)(\boldsymbol{b}^{\top}\boldsymbol{\alpha} - \alpha_t)\\
         && + \delta(\boldsymbol{a}^{\top}\boldsymbol{\alpha} - \alpha_s)(\boldsymbol{b}^{\top}\boldsymbol{1} - 1) + \boldsymbol{a}^{\top}\boldsymbol{\Sigma}\boldsymbol{b} - \boldsymbol{a}^{\top}\boldsymbol{\omega}_t - \boldsymbol{b}^{\top}\boldsymbol{\omega}_s \\
         &&+ \omega_{st}  + (\boldsymbol{a}^{\top}\boldsymbol{\alpha} - \alpha_s)(\boldsymbol{b}^{\top}\boldsymbol{\alpha} - \alpha_t)\bigg].
\end{eqnarray*}The determinant of the MSPE matrix is then $W = W_1W_2-W^2_3$ which needs to be minimized with respect to $\boldsymbol{a}$ and $\boldsymbol{b}$. Upon differentiating $W$ with respect to $ \boldsymbol{a}$ and $ \boldsymbol{b}$ and equating them to vector $ \boldsymbol{0}$ of size $r\times 1$, we obtain 
\begin{equation}\label{a}
(\boldsymbol{\Gamma}\boldsymbol{a} - \boldsymbol{\Delta}_{s})W_2 - 	 (\boldsymbol{\Gamma}\boldsymbol{b} - \boldsymbol{\Delta}_{t})W_3 = \boldsymbol{0}
\end{equation}and
\begin{equation}\label{b}
		(\boldsymbol{\Gamma}\boldsymbol{b} - \boldsymbol{\Delta}_{t})W_1 - 	(\boldsymbol{\Gamma}\boldsymbol{a} - \boldsymbol{\Delta}_{s})W_3 = \boldsymbol{0}, 
	\end{equation}where $\boldsymbol{\Gamma}$, $\boldsymbol{\Delta}_{s}$ and $\boldsymbol{\Delta}_{t}$ are as given in the statement of the theorem. Now, from (\ref{a}) and (\ref{b}), we obtain
\begin{equation}\label{W}
	(\boldsymbol{\Gamma}\boldsymbol{a} - \boldsymbol{\Delta}_{s})(W_1W_2 -W^2_3) = \boldsymbol{0}.
\end{equation}Assuming $W_1W_2 -W^2_3\neq 0$, we readily obtain
$\boldsymbol{a} = \boldsymbol{\Gamma}^{-1} \boldsymbol{\Delta}_{s}$. Similarly, from (\ref{a}) and (\ref{b}), we also obtain $\boldsymbol{b} = \boldsymbol{\Gamma}^{-1} \boldsymbol{\Delta}_{t}$. Hence, the theorem.

\subsection{Interpretation of the quantity $\delta$}
\begin{figure}[h]
	\centering
	\includegraphics[width=10cm,height=12cm,angle=-90,scale=1.15]{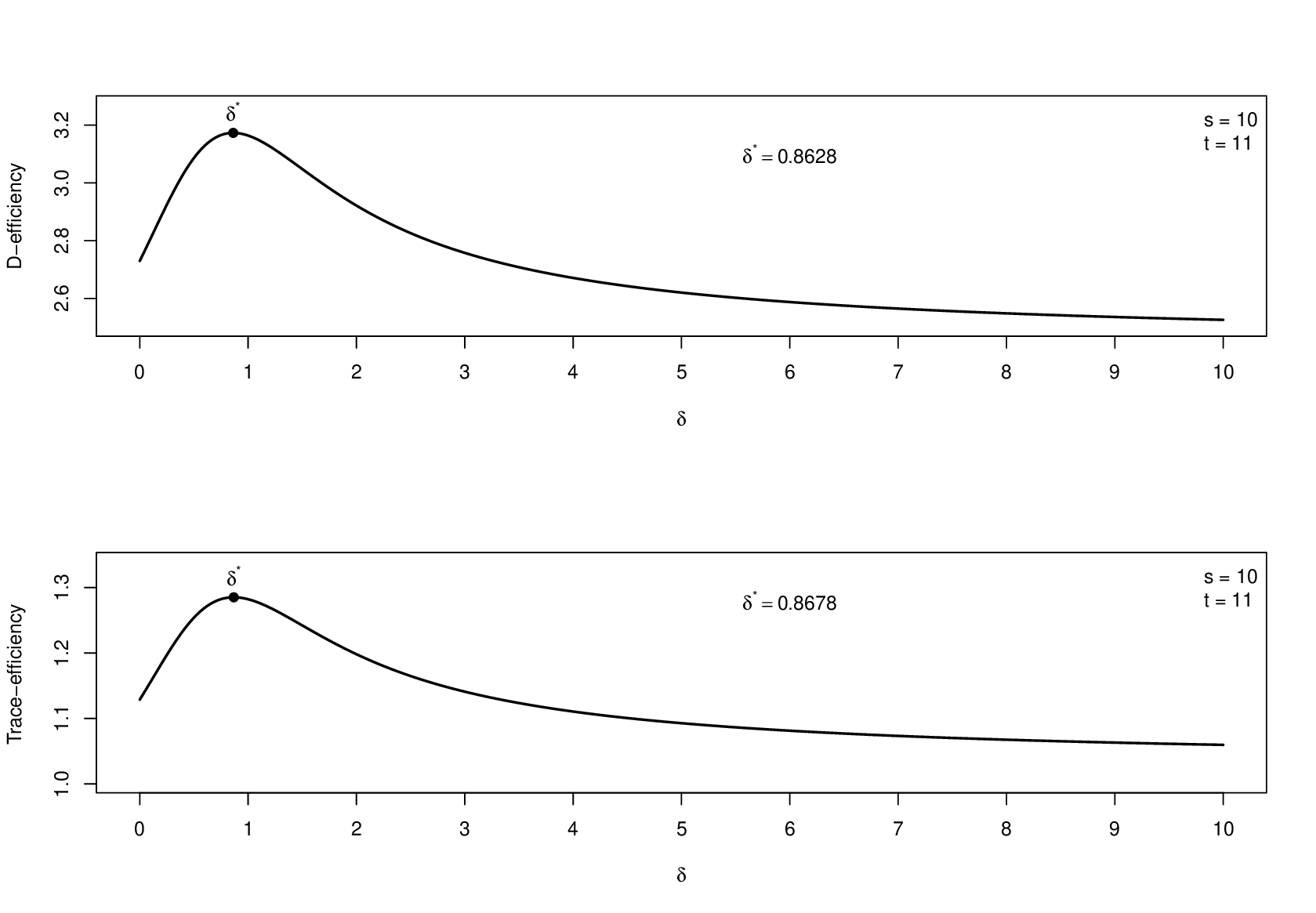}
	\caption{Plot of relative efficiencies of the joint BLIP against joint BLUP for the $s$-th and $t$-th order statistics based on $r=9$ and $n=15$.  }
	\label{Eff_plot_2}
\end{figure}

\begin{figure}[h]
	\centering
	\includegraphics[width=10cm,height=12cm,angle=-90,scale=1.15]{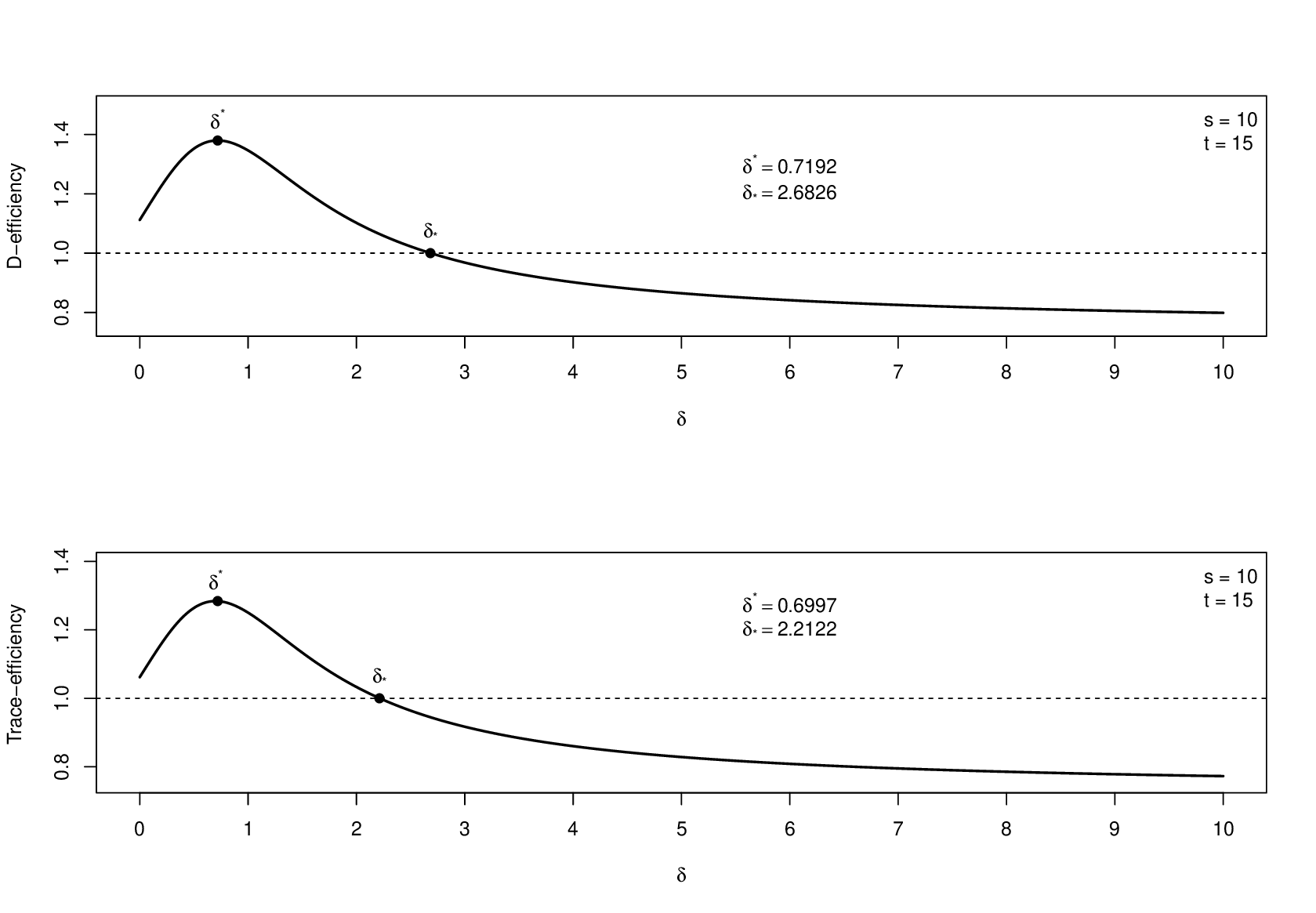}
	\caption{Plot of relative efficiencies of the joint BLIP against joint BLUP for the $s$-th and $t$-th order statistics based on $r=9$ and $n=15$. }
	\label{Eff_plot_3}
\end{figure}

\begin{figure}[h]
	\centering
	\includegraphics[width=10cm,height=12cm,angle=-90,scale=1.15]{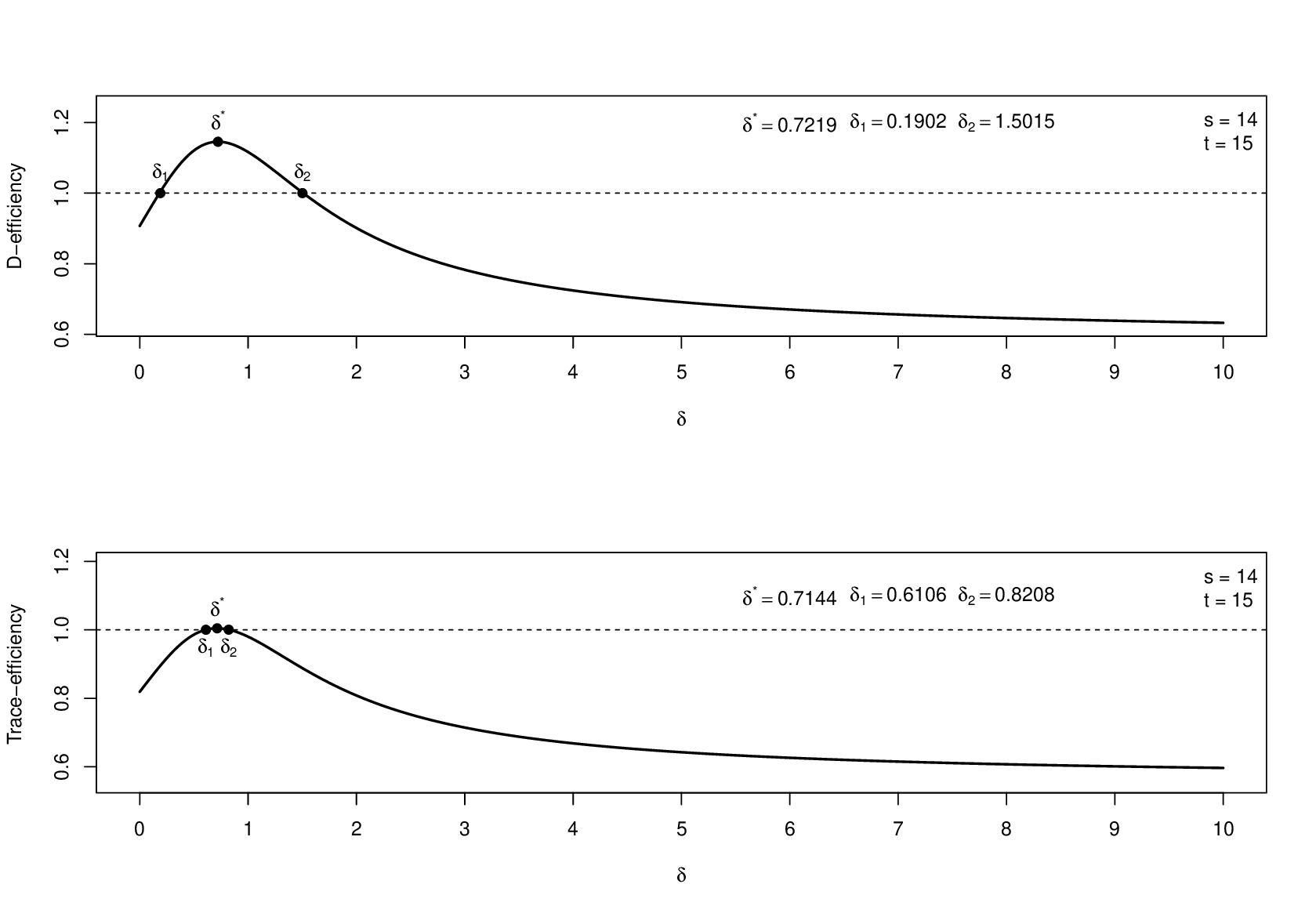}
	\caption{Plot of relative efficiencies of the joint BLIP against joint BLUP for the $s$-th and $t$-th order statistics based on $r=9$ and $n=15$.  }
	\label{Eff_plot_4}
\end{figure}

\paragraph{}
In the case of joint prediction, we define two types of relative efficiencies based on the MSPE matrix as follows:
\begin{eqnarray*}
	\mbox{D-efficiency} &=& \frac{ 
		\mbox{Determinant of the MSPE matrix of $(\hat{X}_{s:n}, \hat{X}_{t:n})$}}{
		\mbox{Determinant of the MSPE matrix of $(\tilde{X}_{s:n}, \hat{X}_{t:n})$}}\\
	&=& \frac{W_{1}W_{2} - W^2_{3}}{W_{11}W_{22} - W^2_{12}}
\end{eqnarray*}and
\begin{eqnarray*}
	\mbox{Trace-efficiency} &=& \frac{ 
		\mbox{Trace of the MSPE matrix of $(\hat{X}_{s:n}, \hat{X}_{t:n})$}}{
		\mbox{Trace of the MSPE matrix of $(\tilde{X}_{s:n}, \hat{X}_{t:n})$}}\\
	&=& \frac{W_{1}+W_{2}}{W_{11}+W_{22} }.
\end{eqnarray*}These relative efficiencies are plotted against $\delta$ for three different pairs of choices: (i) $s=r+1, t=r+2$; (ii) $s=r+1, t=n$; and (iii) $s=n-1, t=n$. Here we took $n=15$ and $s=9$. The corresponding plots are presented in Figures \ref{Eff_plot_2}-\ref{Eff_plot_4}. \\

Both D-efficiency and Trace-efficiency always possess an unique maximum $\delta^*$ which are calculated numerically and indicated in the figures. In Figure \ref{Eff_plot_2}, both efficiencies are always greater than 1 indicating that BLUP performs better than BLIP. In Figure \ref{Eff_plot_3}, BLUP performs better in an interval $(0, \delta_*)$ and then BLIP performs better in the interval $(\delta_*, \infty)$. For this reason, for overall comparative assessment, we define an integrated efficiency measure (IEM) which is simply an average of all efficiencies calculated on a finite interval $(0, \delta_{\tiny{\mbox{max}}})$ for $\delta$. We computed IEM(D-efficiency)  and IEM(Trace-efficiency) for some values of $\delta_{\tiny{\mbox{Max}}}$ and these are presented in Table 2. Table 2 shows both IEM(D-efficiency)  and IEM(Trace-efficiency) are less than 1 which indicates that BLIP has overall better performance in the specified range of $\delta$. In Figure \ref{Eff_plot_4}, it is seen that BLUP performs better in an interval $(\delta_1, \delta_2)$ and BLIP performs better outside that interval. Although the corresponding numerical results are not presented, IEM values are found to be less than 1 in this case as well. Unlike the marginal predictor case in Section 3.1, BLIP is not uniformly better than BLUP in the joint predictor case; but, a practical data-driven guideline on choosing between BLIP or BLUP can be made by estimating $\delta$ using BLUEs of $\mu$ and $\sigma$. Comparing the estimated value of $\delta$, say $\hat{\delta}$, with $\delta_1$ and $\delta_2$, one can compute the gain in efficiency while  using BLIP or BLUP. 
\begin{table}[h]
\centering
\caption{Integrated efficiency measures for D-efficiency and Trace-efficiency for the joint prediction of $X_{10:15}$ and $X_{15:15}$.}
\begin{tabular}{lcccc}\toprule
	$\delta_{\tiny{\mbox{max}}}$ && IEM(D-efficiency) && IEM(Trace-efficiency) \\\bottomrule
	10 && 0.9484 && 0.9024 \\
	50 && 0.8029 && 0.7769 \\
	1000 && 0.7513 && 0.7330\\
	10000 && 0.7486 && 0.7299\\\bottomrule
\end{tabular}	
\end{table}

\subsection{An illustrative example}
\paragraph{}
Let us consider the same data set as presented in Section 3.2. Based on the estimated value $\delta^* = 1.257$,  D- and Trace-efficiencies are computed for three different pairs of choices: (i) $s=10, t=11$; (ii) $s=10, t=15$; and (iii) $s=14, t=15$. The results are summarized in Table 3. Given the data, estimate of $\delta$ indicates that the joint BLUPs are better than joint BLIPs most of the cases except when trace minimizing criterion is opted for determining joint predictors of two extreme order statistics. 
\begin{table}[h]
	\centering
	\scriptsize
	\caption{Summary of D- and Trace-efficiencies based on the data presented by \cite{Bhaumik_2004} and \cite{Krishnamoorthy_2018}. Here, $(\mu^*, \sigma^*) = (2.253, 1.696)$ and $\hat{\delta}=1.328$. }
	\begin{tabular}{lllllllll}\toprule
\multirow{2}{*}{$(s, t)$} && \multicolumn{3}{c}{Determinant based criterion} && \multicolumn{3}{c}{Trace based criterion}\\
 && MSPE(BLIP) & MSPE(BLUP) & D-efficiency && MSPE(BLIP) & MSPE(BLUP) & Trace-efficiency\\\toprule
 (10, 11) && 0.00091 & 0.00029 & 3.137 && 0.0923 & 0.0734 & 1.257 \\
 (10, 15) && 0.0126 & 0.0097 & 1.298 && 0.5324 & 0.4533 & 1.174\\
 (14, 15) && 0.0532 & 0.0510 & 1.043 && 0.7735 & 0.8399 & 0.9209\\\bottomrule
	\end{tabular}
\end{table}

\subsection{Complete MSPE matrix dominance of BLIPs}
\paragraph{}
Let us consider the problem of predicting the random quantity $Y = lX_{s:n} + kX_{t:n},$ which is a linear combination of two future order statistics $X_{s:n}$ and $X_{t:n}$, with $l$ and $k$ being two arbitrary fixed constants. Let us then assume a linear predictor for $Y$ as $$\hat{Y}=\boldsymbol{c}^{\top}\boldsymbol{X},$$where the coefficient vector $\boldsymbol{c}=(c_1, \ldots, c_r)^{\top}_{r\times 1}$ needs to be suitably determined by minimizing the mean squared error of the predictor $\hat{Y}$. The mean squared predictive error of $\hat{Y}$ is given by 
\begin{eqnarray*}
	W &=& \sigma^2 \bigg[\delta^2(\boldsymbol{c}^{\top}\boldsymbol{1}-l-k)^2 +2\delta(\boldsymbol{c}^{\top}\boldsymbol{1}-l-k)(\boldsymbol{c}^{\top}\boldsymbol{\alpha}-l\alpha_s-k\alpha_t)\\
	&&+\boldsymbol{c}^{\top}\boldsymbol{\Sigma}\boldsymbol{c} + (\boldsymbol{c}^{\top}\boldsymbol{\alpha}-l\alpha_s-k\alpha_t)^2 - 2l\boldsymbol{c}^{\top}\boldsymbol{\omega}_s-2k\boldsymbol{c}^{\top}\boldsymbol{\omega}_s-2kl\omega_{st}\\
	&& +l^2\omega_{ss} + k^2\omega_{tt}\bigg].
\end{eqnarray*}Now, taking derivative of $W$ with respect to $\boldsymbol{c}$ and then equating it to null vector $\boldsymbol{0}$, we obtain
\begin{eqnarray*}
	&&\delta^2(\boldsymbol{c}^{\top}\boldsymbol{1}-l-k)\boldsymbol{1}  + \delta \boldsymbol{1} (\boldsymbol{c}^{\top}\boldsymbol{\alpha}-l\alpha_s-k\alpha_t)+ \delta (\boldsymbol{c}^{\top}\boldsymbol{1}-l-k)\boldsymbol{\alpha}\\
&&+	\boldsymbol{\Sigma}\boldsymbol{c} + (\boldsymbol{c}^{\top}\boldsymbol{\alpha}-l\alpha_s-k\alpha_t)\boldsymbol{\alpha} - l\boldsymbol{\omega}_s-k\boldsymbol{\omega}_t = \boldsymbol{0},
	\end{eqnarray*}which yields
\begin{equation*}
	\boldsymbol{\Gamma}\boldsymbol{c} = l \boldsymbol{\Delta}_s + k\boldsymbol{\Delta}_t,
\end{equation*}where $\boldsymbol{\Gamma}$, $\boldsymbol{\Delta}_s$ and $\boldsymbol{\Delta}_t$ are as defined in Theorem 2. As $\boldsymbol{\Gamma}$ is invertible, we then obtain
\begin{equation*}
	\boldsymbol{c} = l \boldsymbol{\Gamma}^{-1}\boldsymbol{\Delta}_s + k\boldsymbol{\Gamma}^{-1}\boldsymbol{\Delta}_t,
\end{equation*}and consequently, 
\begin{equation*}
	\hat{Y} = l \hat{X}_{r:n} + k\hat{X}_{t:n},
\end{equation*}where $\hat{X}_{s:n}$ and $\hat{X}_{t:n}$ are the joint BLIPs of $X_{s:n}$ and $X_{t:n}$, respectively, based on $\boldsymbol{X}$. As a result, we have
\begin{equation*}
	\mbox{Var}( l \hat{X}_{s:n} + k\hat{X}_{t:n})\leq \mbox{Var}( l X^{*}_{s:n} + kX^{*}_{t:n})
\end{equation*}for any other joint linear invariant predictors $X^{*}_{s:n}$ and $X^{*}_{t:n}$ of $X_{s:n}$ and $X_{t:n}$. This readily implies
\begin{eqnarray*}
	&&\begin{bmatrix} l &  k\end{bmatrix}\begin{bmatrix} \mbox{MVar}(\hat{X}_{s:n}) &  \mbox{MCov}(\hat{X}_{s:n}, \hat{X}_{t:n})\\ \mbox{MCov}(\hat{X}_{s:n}, \hat{X}_{t:n}) & \mbox{MVar}(\hat{X}_{t:n})\end{bmatrix}\begin{bmatrix} l \\ k\end{bmatrix}\\
	&&\leq \begin{bmatrix} l &  k\end{bmatrix}\begin{bmatrix} \mbox{MVar}(X^{*}_{s:n}) &  \mbox{MCov}(X^{*}_{s:n}, X^{*}_{t:n})\\ \mbox{MCov}(X^{*}_{s:n}, X^{*}_{t:n}) & \mbox{MVar}(X^{*}_{t:n})\end{bmatrix}\begin{bmatrix} l \\ k\end{bmatrix},
\end{eqnarray*}where MVar and MCov stand for mean squared predictive error variance and mean squared predictive error covariance of the joint predictors, respectively. This establishes the property that the BLIPs of $X_{s:n}$ and $X_{t:n}$ possess complete MSPE matrix dominance in the class of all linear invariant predictors of $X_{s:n}$ and $X_{t:n}$, which is a more general property than  trace-efficiency and determinant-efficiency.

\section{Extension to prediction of $\ell$ order statistics}
\paragraph{}
Let us now consider the BLIPs of any $\ell$ future order statistics $(X_{s_1:n}, X_{s_2:n},$ $\ldots, X_{s_{\ell}:n}),$ for $r<s_1<s_2<\ldots<s_{\ell}\leq n,$ simultaneously. We then have the following general result.
\begin{thm}
	The simultaneous best linear invariant predictors of any $\ell$ future order statistics are identical to their corresponding marginal predictors. 
\end{thm}
\noindent\textbf{Proof:} Let us assume that the BLIPs of $\ell$ future order statistics are of the form 
\begin{equation}\label{BLUPai}
	\hat{X}_{s_i:n} = \boldsymbol{a}^{\top}_i \boldsymbol{X}, ~i \in \{1, 2,\ldots, \ell\},
\end{equation}where $\boldsymbol{a}^{,}_i$s are the coefficient vectors of size $r\times1$ that need to be suitably determined.  The corresponding mean squared predictive error matrix is then $$\boldsymbol{W} =\big(\big(W_{ij}\big)\big)_{i, j=1}^{\ell},$$where
\begin{eqnarray*}
	W_{ii} &=& \sigma^2\bigg[  \delta^2(\boldsymbol{a}^{\top}_i\boldsymbol{1} - 1)^2 +2\delta (\boldsymbol{a}^{\top}_i\boldsymbol{1} - 1)(\boldsymbol{a}^{\top}_i\boldsymbol{\alpha} - \alpha_{s_i})\\
	&& ~~~~~~~~~~+\boldsymbol{a}^{\top}_i\boldsymbol{\Sigma}\boldsymbol{a}_i -2\boldsymbol{a}^{\top}_i\boldsymbol{\omega}_{s_i} - \omega_{s_is_i}  +  (\boldsymbol{a}^{\top}_i\boldsymbol{\alpha} - \alpha_{s_i})^2\bigg],\quad i\in\{1, 2,\ldots, \ell\},
\end{eqnarray*}and
\begin{eqnarray*}
	W_{ij}  &=& \sigma^2\bigg[   \delta^2(\boldsymbol{a}^{\top}_i\boldsymbol{1} - 1)(\boldsymbol{a}^{\top}_j\boldsymbol{1} - 1)+\delta(\boldsymbol{a}^{\top}_i\boldsymbol{1} - 1)(\boldsymbol{a}^{\top}_j\boldsymbol{\alpha} - \alpha_{s_j})\\
	&& + \delta(\boldsymbol{a}^{\top}_i\boldsymbol{\alpha} - \alpha_{s_i})(\boldsymbol{a}^{\top}_j\boldsymbol{1} - 1) + \boldsymbol{a}^{\top}_i\boldsymbol{\Sigma}\boldsymbol{a}_j - \boldsymbol{a}^{\top}_i\boldsymbol{\omega}_{s_j} - \boldsymbol{a}^{\top}_j\boldsymbol{\omega}_{s_i} \\
	&&+ \omega_{s_is_j}  + (\boldsymbol{a}^{\top}_i\boldsymbol{\alpha} - \alpha_{s_i})(\boldsymbol{a}^{\top}_j\boldsymbol{\alpha} - \alpha_{s_j})\bigg], \quad 1\leq i < j\leq \ell.
\end{eqnarray*}We observe that $\boldsymbol{W}$ is symmetric, i.e., $W_{ij} = W_{ji}$, and further that, each coefficient vector $\boldsymbol{a}_i$ appears in only one row and one column. For instance, $\boldsymbol{a}_i$ appears only in the $i$th row and the $i$th column. Let us further denote
\begin{equation*}
	\frac{\partial}{\partial \boldsymbol{a}_i } W_{ii} = 2(\boldsymbol{\Gamma}\boldsymbol{a}_i - \boldsymbol{\Delta}_i), \quad i\in\{1, 2,\ldots, \ell\},
\end{equation*}
and
\begin{eqnarray*}
	\frac{\partial}{\partial \boldsymbol{a}_i } W_{ij} &=& \boldsymbol{\Gamma}\boldsymbol{a}_j - \boldsymbol{\Delta}_j,\\
	\frac{\partial}{\partial \boldsymbol{a}_j} W_{ij} &=& \boldsymbol{\Gamma}\boldsymbol{a}_i - \boldsymbol{\Delta}_i, \quad 1\leq i \neq j\leq \ell,
\end{eqnarray*}where $\boldsymbol{\Gamma} = \boldsymbol{\Sigma} + (\boldsymbol{\alpha}+\delta\boldsymbol{1})  (\boldsymbol{\alpha}+\delta\boldsymbol{1})^{\top}$ and $\boldsymbol{\Delta}_{s} = \boldsymbol{\omega}_s + (\alpha_s + \delta)(\boldsymbol{\alpha} + \delta\boldsymbol{1})$.\\

In addition, let us use $|\boldsymbol{W}|$ to denote the determinant of $\boldsymbol{W}$ which needs to be minimized with respect to $\boldsymbol{a}_i$, $i = 1,\cdots, \ell$. Taking derivative of 
$|\boldsymbol{W}|$ with respect to $\boldsymbol{a}_1$, for example, we obtain
\begin{eqnarray*}
		\frac{\partial}{\partial \boldsymbol{a}_1 }|\boldsymbol{W}|  &=& \sigma^2\frac{\partial}{\partial \boldsymbol{a}_1 } \begin{vmatrix}
		W_{11} & W_{12} & \ldots & W_{1\ell}\\
		W_{12} & W_{22} & \ldots & W_{2\ell}\\
		\vdots & \vdots &  & \vdots \\
		W_{1\ell} & W_{2\ell} & \ldots & W_{\ell\ell}\\
	\end{vmatrix}\\
&=& \sigma^2\begin{vmatrix}
	\frac{1}{2}\frac{\partial W_{11}}{\partial \boldsymbol{a}_1 } & 	\frac{\partial W_{12}}{\partial \boldsymbol{a}_1 } & \ldots & 	\frac{\partial W_{1\ell}}{\partial \boldsymbol{a}_1 }\\
	W_{12} & W_{22} & \ldots & W_{2\ell}\\
	\vdots & \vdots &  & \vdots \\
	W_{1\ell} & W_{2\ell} & \ldots & W_{\ell\ell}\\
\end{vmatrix}+
\sigma^2\begin{vmatrix}
	\frac{1}{2}\frac{\partial W_{11}}{\partial \boldsymbol{a}_1 } & W_{12} & \ldots & 	W_{1\ell}\\
	\frac{\partial W_{12}}{\partial \boldsymbol{a}_1 } & W_{22} & \ldots & W_{2\ell}\\
	\vdots & \vdots &  & \vdots \\
	\frac{\partial W_{1\ell}}{\partial \boldsymbol{a}_1 } & W_{2\ell} & \ldots & W_{\ell\ell}\\
\end{vmatrix} \\
&=& 2\sigma^2\begin{vmatrix}
	\frac{1}{2}\frac{\partial W_{11}}{\partial \boldsymbol{a}_1 } & 	\frac{\partial W_{12}}{\partial \boldsymbol{a}_1 } & \ldots & 	\frac{\partial W_{1\ell}}{\partial \boldsymbol{a}_1 }\\
	W_{12} & W_{22} & \ldots & W_{2\ell}\\
	\vdots & \vdots &  & \vdots \\
	W_{1\ell} & W_{2\ell} & \ldots & W_{\ell\ell}\\
\end{vmatrix}\quad\mbox{(due to symmetry of the determinants)}
\end{eqnarray*}
\begin{eqnarray*}
&=& 2\sigma^2\begin{vmatrix}
\boldsymbol{\Gamma}\boldsymbol{a}_1 - \boldsymbol{\Delta}_1 & 	\boldsymbol{\Gamma}\boldsymbol{a}_2 - \boldsymbol{\Delta}_2  & \ldots & 	\boldsymbol{\Gamma}\boldsymbol{a}_{\ell} - \boldsymbol{\Delta}_{\ell} \\
	W_{12} & W_{22} & \ldots & W_{2\ell}\\
	\vdots & \vdots &  & \vdots \\
	W_{1\ell} & W_{2\ell} & \ldots & W_{\ell\ell}\\
\end{vmatrix}\\
&=& 2\sigma^2\begin{vmatrix}
	\boldsymbol{\Gamma}\boldsymbol{a}_1  & 	\boldsymbol{\Gamma}\boldsymbol{a}_2  & \ldots & 	\boldsymbol{\Gamma}\boldsymbol{a}_{\ell} \\
	W_{12} & W_{22} & \ldots & W_{2\ell}\\
	\vdots & \vdots &  & \vdots \\
	W_{1\ell} & W_{2\ell} & \ldots & W_{\ell\ell}\\
\end{vmatrix} - 
2\sigma^2\begin{vmatrix}
	 \boldsymbol{\Delta}_1 &  \boldsymbol{\Delta}_2  & \ldots & 	\boldsymbol{\Delta}_{\ell} \\
	W_{12} & W_{22} & \ldots & W_{2\ell}\\
	\vdots & \vdots &  & \vdots \\
	W_{1\ell} & W_{2\ell} & \ldots & W_{\ell\ell}\\
\end{vmatrix}\\
&=& 2\sigma^2 M_1 - 2\sigma^2 M_2,
\end{eqnarray*}where
\begin{equation*}
	M_1 = \begin{vmatrix}
		\boldsymbol{\Gamma}\boldsymbol{a}_1  & 	\boldsymbol{\Gamma}\boldsymbol{a}_2  & \ldots & 	\boldsymbol{\Gamma}\boldsymbol{a}_{\ell} \\
		W_{12} & W_{22} & \ldots & W_{2\ell}\\
		\vdots & \vdots &  & \vdots \\
		W_{1\ell} & W_{2\ell} & \ldots & W_{\ell\ell}\\
	\end{vmatrix}\mbox{~and~}
M_2 = \begin{vmatrix}
	\boldsymbol{\Delta}_1 &  \boldsymbol{\Delta}_2  & \ldots & 	\boldsymbol{\Delta}_{\ell} \\
	W_{12} & W_{22} & \ldots & W_{2\ell}\\
	\vdots & \vdots &  & \vdots \\
	W_{1\ell} & W_{2\ell} & \ldots & W_{\ell\ell}\\
\end{vmatrix}.
\end{equation*}Now, expanding the determinant $M_1$ by its first row, we readily find  it to be 
\begin{equation*}
	\boldsymbol{\Gamma}\boldsymbol{a}_1 C^{11} + \boldsymbol{\Gamma}\boldsymbol{a}_2 C^{12}+\cdots+ \boldsymbol{\Gamma}\boldsymbol{a}_{\ell} C^{1\ell},
\end{equation*}where $C^{ij}$ is the co-factor of $W_{ij}$. Similarly, expanding the determinant $M_2$ by its first row, we obtain it to be
\begin{equation*}
	\boldsymbol{\Delta}_1  C^{11} + \boldsymbol{\Delta}_2 C^{12}+\cdots+\boldsymbol{\Delta}_{\ell} C^{1\ell}.
\end{equation*}Therefore, $	\frac{\partial}{\partial \boldsymbol{a}_1 }|\boldsymbol{W}| = 0$ yields the following normal equation
\begin{equation*}
	\begin{bmatrix}
	C^{11} & C^{12} &\cdots& C^{1\ell}
  \end{bmatrix}
\begin{bmatrix}
	\boldsymbol{\Gamma}\boldsymbol{a}_1 \\ \boldsymbol{\Gamma}\boldsymbol{a}_2 \\\vdots \\ \boldsymbol{\Gamma}\boldsymbol{a}_{\ell} 
\end{bmatrix} =
\begin{bmatrix}
	C^{11} & C^{12} &\cdots& C^{1\ell}
\end{bmatrix}
\begin{bmatrix}
\boldsymbol{\Delta}_1 \\ \boldsymbol{\Delta}_2 \\\vdots \\ \boldsymbol{\Delta}_{\ell}
\end{bmatrix}.	
\end{equation*}
Similarly, by taking derivative of  $|\boldsymbol{W}|$ with respect to $\boldsymbol{a}_i$, $i = 2,\cdots, \ell$, and proceeding exactly as above, we can generate the following system of equations 
\begin{eqnarray*}
		\begin{bmatrix}
		C^{21} & C^{22} &\cdots& C^{2\ell}
	\end{bmatrix}
	\begin{bmatrix}
		\boldsymbol{\Gamma}\boldsymbol{a}_1 \\ \boldsymbol{\Gamma}\boldsymbol{a}_2 \\\vdots \\ \boldsymbol{\Gamma}\boldsymbol{a}_{\ell} 
	\end{bmatrix}& =&
	\begin{bmatrix}
		C^{21} & C^{22} &\cdots& C^{2\ell}
	\end{bmatrix}
	\begin{bmatrix}
		\boldsymbol{\Delta}_1 \\ \boldsymbol{\Delta}_2 \\\vdots \\ \boldsymbol{\Delta}_{\ell}
	\end{bmatrix},\\
\vdots~~~~~~~~~~~~~~~~~~~~~~&=&~~~~~~~~~~~~~~~~~\vdots\\
	\begin{bmatrix}
	C^{\ell1} & C^{\ell2} &\cdots& C^{\ell\ell}
\end{bmatrix}
\begin{bmatrix}
	\boldsymbol{\Gamma}\boldsymbol{a}_1 \\ \boldsymbol{\Gamma}\boldsymbol{a}_2 \\\vdots \\ \boldsymbol{\Gamma}\boldsymbol{a}_{\ell} 
\end{bmatrix} &=&
\begin{bmatrix}
	C^{\ell1} & C^{\ell2} &\cdots& C^{\ell\ell}
\end{bmatrix}
\begin{bmatrix}
	\boldsymbol{\Delta}_1 \\ \boldsymbol{\Delta}_2 \\\vdots \\ \boldsymbol{\Delta}_{\ell}
\end{bmatrix}.	
\end{eqnarray*}Then, all these $\ell$ equations can be written in a combined form as 
\begin{equation}\label{Cmatrix}
		\begin{bmatrix}
		C^{11} & C^{12} & \ldots & C^{1\ell}\\
		C^{12} & C^{22} & \ldots & C^{2\ell}\\
		\vdots & \vdots &  & \vdots \\
		C^{1\ell} & C^{2\ell} & \ldots & C^{\ell\ell}\\
	\end{bmatrix}\begin{bmatrix}
	\boldsymbol{\Gamma}\boldsymbol{a}_1 \\ \boldsymbol{\Gamma}\boldsymbol{a}_2 \\\vdots \\ \boldsymbol{\Gamma}\boldsymbol{a}_{\ell} 
\end{bmatrix}
=
	\begin{bmatrix}
	C^{11} & C^{12} & \ldots & C^{1\ell}\\
	C^{12} & C^{22} & \ldots & C^{2\ell}\\
	\vdots & \vdots &  & \vdots \\
	C^{1\ell} & C^{2\ell} & \ldots & C^{\ell\ell}\\
\end{bmatrix}\begin{bmatrix}
\boldsymbol{\Delta}_1 \\ \boldsymbol{\Delta}_2 \\\vdots \\ \boldsymbol{\Delta}_{\ell}
\end{bmatrix}.
\end{equation}With $\boldsymbol{C} =  \big(\big(C^{ij}\big)\big)_{i, j=1}^l$ denoting the adjoint matrix of $\boldsymbol{W}$, it is known that $\boldsymbol{C}=|\boldsymbol{W}|\boldsymbol{W}^{-1}.$ As $\boldsymbol{W}$ is positive-definite and is invertible, so is $\boldsymbol{C}$. Thence, by pre-multiplying (\ref{Cmatrix}) by $\boldsymbol{C}^{-1}$ on both sides, we readily obtain
\begin{equation}\label{gammaeqn}
	\begin{bmatrix}
		\boldsymbol{\Gamma}\boldsymbol{a}_1 \\ \boldsymbol{\Gamma}\boldsymbol{a}_2 \\\vdots \\ \boldsymbol{\Gamma}\boldsymbol{a}_{\ell} 
	\end{bmatrix} = \begin{bmatrix}
	\boldsymbol{\Delta}_1 \\ \boldsymbol{\Delta}_2 \\\vdots \\ \boldsymbol{\Delta}_{\ell}
\end{bmatrix}.
\end{equation}Note that $\boldsymbol{\Gamma}$ is a sum of two positive-definite quadratic forms and is therefore positive-definite. So, its inverse exists uniquely, and hence, the solution of (\ref{gammaeqn}) is simply
\begin{equation*}
	\boldsymbol{a}_i = \boldsymbol{\Gamma}^{-1}\boldsymbol{\Delta}_i,\quad i\in\{1, 2,\ldots, \ell\},
\end{equation*}which completes the proof of the theorem. 

\begin{cor}
	As done in Section 4.2, one can easily establish that the simultaneous BLIPs derived in Theorem 3 also possess the complete MSPE matrix dominance property. 
\end{cor}

\section{BLIPs in scale family of distributions}
\paragraph{}
So far, we have discussed the prediction problem for a general location-scale family of distributions. In a much simpler way, analogous results can be developed for scale family of  distributions. Let us assume that the parent distribution of the first $r$ observed order statistics $\boldsymbol{X}$ belongs to the scale family whose probability density function is given by
\begin{equation*}
	\frac{1}{\sigma} f\left(\frac{x}{\sigma} \right),\quad \sigma>0.
\end{equation*}Now, let us denote $\alpha_i = \mbox{E}[Z_{i:n}] = \mbox{E}[X_{i:n}/\sigma]$ for $i\in\{1,\ldots, r\}$, $\boldsymbol{\alpha}=(\alpha_1, \ldots, \alpha_r)^{\top}$ and $\boldsymbol{\Sigma}$ as the $r\times r$ covariance matrix of $\boldsymbol{Z} = (Z_{1:n},\ldots, Z_{r:n})^{\top},$ assumed to be positive definite. Then, the marginal and simultaneous BLIPs are as presented in the following results.
\begin{thm}
	The marginal best linear invariant predictor $\tilde{X}_{s:n}$ of $X_{s:n}$, determined by minimizing the mean squared predictive error of $\tilde{X}_{s:n}$, is of the form $\tilde{X}_{s:n} = \boldsymbol{a}^{\top}\boldsymbol{X}$ in which the coefficient vector $\boldsymbol{a}=(a_1, \ldots, a_r)^{\top}_{r\times 1}$ is given by
	\begin{equation*}
		\boldsymbol{a} = (\boldsymbol{\Sigma} + \boldsymbol{\alpha}\boldsymbol{\alpha}^{\top})^{-1} (\boldsymbol{\omega}_s + \alpha_s\boldsymbol{\alpha}).
	\end{equation*}
\end{thm}
\noindent\textbf{Proof:} Similar to Theorem 1 for the location-scale family, the proof follows in this case by minimizing the mean squared predictive error given by
\begin{equation*}
\sigma^2\bigg[	\boldsymbol{a}^{\top}\boldsymbol{\Sigma}\boldsymbol{a} -2\boldsymbol{a}^{\top}\boldsymbol{\omega}_s + \omega_{ss} + (\boldsymbol{a}^{\top}\boldsymbol{\alpha}- \alpha_s)^2\bigg].
\end{equation*}

\begin{thm}
	The simultaneous best linear invariant predictors of $\ell$ future order statistics are identical to their corresponding marginal predictors. 
\end{thm}
\noindent\textbf{Proof:} The proof follows exactly in the same way as in the proof of Theorem 3 for the location-scale family by considering the mean squared predictive error matrix in this case as
\begin{equation*}
	\begin{bmatrix} W_{1} &  W_{3}\\ W_{3} & W_{2}\end{bmatrix},
\end{equation*}where 
\begin{eqnarray*}
	W_1 &=& \sigma^2\bigg[  \boldsymbol{a}^{\top}\boldsymbol{\Sigma}\boldsymbol{a} -2\boldsymbol{a}^{\top}\boldsymbol{\omega}_s + \omega_{ss}  +  (\boldsymbol{a}^{\top}\boldsymbol{\alpha} - \alpha_s)^2\bigg],\\
	W_2 &=&  \sigma^2\bigg[ \boldsymbol{b}^{\top}\boldsymbol{\Sigma}\boldsymbol{b} -2\boldsymbol{b}^{\top}\boldsymbol{\omega}_t + \omega_{tt}  +  (\boldsymbol{b}^{\top}\boldsymbol{\alpha} - \alpha_t)^2\bigg],\\
	W_3 &=& \sigma^2\bigg[ \boldsymbol{a}^{\top}\boldsymbol{\Sigma}\boldsymbol{b} - \boldsymbol{a}^{\top}\boldsymbol{\omega}_t - \boldsymbol{b}^{\top}\boldsymbol{\omega}_s+ \omega_{st}  + (\boldsymbol{a}^{\top}\boldsymbol{\alpha} - \alpha_s)(\boldsymbol{b}^{\top}\boldsymbol{\alpha} - \alpha_t)\bigg].
\end{eqnarray*}

\begin{cor}
	All the associated properties of BLIPs presented earlier for the location-scale family in Section 4 can be shown to hold here for the scale family as well. 
\end{cor}

\section{Concluding remarks}
\paragraph{}
In this article, we have presented explicit expressions for the simultaneous BLIPs of any $\ell$ future order statistics and have established that the simultaneous BLIPs are the same as the marginal BLIPs. The advantage of using BLIP over BLUP in marginal prediction case has been demonstrated. Moreover, in the simultaneous prediction case, a practical data-driven approach for choosing between BLIP and BLUP has been discussed as well.

\section*{Acknowledgments}
\paragraph{}
The first author thanks the Natural Sciences and Engineering Research Council of Canada for funding this research through an Individual Discovery Grant (RGPIN-2020-0633).

\bibliographystyle{apalike}
\bibliography{ritwik_ref}

\end{document}